\documentclass[a4paper,12pt]{amsart}
\usepackage{amssymb}
\usepackage{ifthen}
\usepackage[usenames]{color}
\usepackage{graphicx}
\nonstopmode \numberwithin{equation}{section}
\newtheorem{thm}{Theorem}[section]

\newtheorem{cor}{Corollary}[section]
\newtheorem{lem}{Lemma}[section]
\newtheorem{prop}{Proposition}[section]

\newtheorem{claim}{Claim}[section]
\newtheorem{conj}[equation]{Conjecture}

\theoremstyle{definition}
\newtheorem{defn}{Definition}[section]
\newtheorem{examp}{Example}[section]
\newtheorem{prob}[equation]{Problem}
\newtheorem{ques}[equation]{Question}
\newtheorem{rem}{Remark}[section]

\newtheorem{case}{Case}[section]
\newtheorem{step}{Step}[section]
\newtheorem*{mysolution}{Solution}

\newcounter {own}
\def\theown {\thesection       .\arabic{own}}

\newenvironment{pf}[1][]{%
 \vskip 3mm
 \noindent
 \ifthenelse{\equal{#1}{}}%
  {{\slshape Proof. }}%
  {{\slshape #1.} }%
 }%
{\qed\bigskip}

\newcounter{alphabet}
\newcounter{tmp}
\newenvironment{Thm}[1][]{\refstepcounter{alphabet}%
\bigskip%
\noindent%
{\bf Theorem \Alph{alphabet}}%
\ifthenelse{\equal{#1}{}}{}{ (#1)}%
{\bf .} \itshape}{\vskip 8pt}

\makeatletter
\newcommand{\Ref}[1]{\@ifundefined{r@#1}{}{\setcounter{tmp}{\ref{#1}}\Alph{tmp}}}
\makeatother

\newcommand{\IR}{{\mathbb R}}

\newcommand{\IC}{{\mathbb C}}

\newcommand{\ID}{{\mathbb D}}
\newcommand{\IB}{{\mathbb B}}

\def\be{\begin{equation}}
\def\ee{\end{equation}}

\newcommand{\ben}{\begin{enumerate}}
\newcommand{\een}{\end{enumerate}}

\newcommand{\blem}{\begin{lem}}
\newcommand{\elem}{\end{lem}}
\newcommand{\bthm}{\begin{thm}}
\newcommand{\ethm}{\end{thm}}
\newcommand{\bcor}{\begin{cor}}
\newcommand{\ecor}{\end{cor}}
\newcommand{\beg}{\begin{examp}}
\newcommand{\eeg}{\end{examp}}
\newcommand{\begs}{\begin{examples}}
\newcommand{\eegs}{\end{examples}}
\newcommand{\bdefe}{\begin{defn}}
\newcommand{\edefe}{\end{defn}}

\newcommand{\bca}{\begin{case}}
\newcommand{\eca}{\end{case}}

\newcommand{\bprob}{\begin{prob}}
\newcommand{\eprob}{\end{prob}}
\newcommand{\bques}{\begin{ques}}
\newcommand{\eques}{\end{ques}}
\newcommand{\bei}{\begin{itemize}}
\newcommand{\eei}{\end{itemize}}
\newcommand{\bcl}{\begin{claim}}
\newcommand{\ecl}{\end{claim}}

\newcommand{\bsol}{\begin{mysolution}}
\newcommand{\esol}{\end{mysolution}}

\newcommand{\bcon}{\begin{conj}}
\newcommand{\econ}{\end{conj}}
\newcommand{\bcons}{\begin{conjs}}
\newcommand{\econs}{\end{conjs}}
\newcommand{\bprop}{\begin{prop}}
\newcommand{\eprop}{\end{prop}}
\newcommand{\br}{\begin{rem}}
\newcommand{\er}{\end{rem}}
\newcommand{\brs}{\begin{rems}}
\newcommand{\ers}{\end{rems}}
\newcommand{\bo}{\begin{obser}}
\newcommand{\eo}{\end{obser}}
\newcommand{\bos}{\begin{obsers}}
\newcommand{\eos}{\end{obsers}}
\newcommand{\bpf}{\begin{pf}}
\newcommand{\epf}{\end{pf}}
\newcommand{\ba}{\begin{array}}
\newcommand{\ea}{\end{array}}
\newcommand{\beq}{\begin{eqnarray}}
\newcommand{\beqq}{\begin{eqnarray*}}
\newcommand{\eeq}{\end{eqnarray}}
\newcommand{\eeqq}{\end{eqnarray*}}

\newcommand{\ds}{\displaystyle}

\newcounter{minutes}
\divide\time by 60
\newcounter{hours}
\multiply\time by 60 \addtocounter{minutes}{-\time}

\begin{document}
\bibliographystyle{amsplain}
\title [Representation formula and bi-Lipschitz continuity of solutions to inhomogeneous biharmonic Dirichlet problems in the unit disk]
{Representation formula and bi-Lipschitz continuity of solutions to inhomogeneous biharmonic Dirichlet problems in the unit disk}

\thanks{
File:~\jobname .tex,
          printed: \number\day-\number\month-\number\year,
          \thehours.\ifnum\theminutes<10{0}\fi\theminutes}
\author{Peijin Li}
\address{P. Li, Department of Mathematics,
Hunan First Normal University, Changsha, Hunan 410205, People's Republic of China}
\email{wokeyi99@163.com}
\author{Saminathan Ponnusamy
}
\address{S. Ponnusamy, Stat-Math Unit,
Indian Statistical Institute (ISI), Chennai Centre,
110, Nelson Manickam Road,
Aminjikarai, Chennai, 600 029, India.}
\email{samy@iitm.ac.in}


\begin{abstract}
The aim of this paper is twofold. First, we establish the representation formula and the uniqueness of the solutions to a class of
inhomogeneous biharmonic Dirichlet problems, and then prove the bi-Lipschitz continuity of the solutions.
\end{abstract}


\maketitle
\pagestyle{myheadings}
\markboth{Peijin Li and Saminathan Ponnusamy
}{Representation formula and bi-Lipschitz continuity of solutions}

\section{Introduction and statement of the main results}\label{sec-1}

Let $\mathbb{C}$ denote the complex plane. For $a\in$ $\mathbb{C}$, let $\mathbb {D}(a,r)=\{z:|z-a|<r\}$, where  $r>0$,
and $\mathbb {D}_{r}=\mathbb {D}(0,r)$. In particular, let $\mathbb {D}=\mathbb {D}_{1}$ and $\mathbb{T}=\partial\mathbb{D}$, the boundary of $\mathbb{D}$.
For domains $D$ and $\Omega$ be domains in $\mathbb{C}$,  a function $p: \,D\to \Omega$ is said to be {\it $L_{1}$-Lipschitz} (resp. {\it $L_2$-co-Lipschitz})
if for all $z,$ $w\in D$,
$$|p(z)-p(w)|\leq L_{1}|z-w|\;\; ( \text{resp.}\;|p(z)-p(w)|\geq L_2|z-w|)
$$
for some positive constants $L_{1}$ and $L_2$.  We say that $p$ is {\it bi-Lipschitz} if it is both Lipschitz and co-Lipschitz.

The main aim of this paper is to discuss the representation formula, the uniqueness and the bi-Lipschitz continuity of the
solutions to the following inhomogeneous biharmonic Dirichlet problem (briefly, IBDP in the following):
\be\label{eq-16}
\begin{cases} {\Delta}^2 \Phi=g &
\mbox{ in }\ds  \ID,\\
\ds \Phi=f & \mbox{ on }\ds \mathbb{T},\\
\ds \partial_n \Phi=h & \mbox{ on }\ds \mathbb{T},
\end{cases}
\ee
where $\Delta$ denotes the Laplacian given by
$$\Delta=\frac{\partial^{2}}{\partial
z\partial\overline{z}}=\frac{1}{4}\left(\frac{\partial ^2}{\partial x^2} +
\frac{\partial ^2}{\partial y^2}\right),
$$
$\partial_n$ denotes the differentiation in the inward normal direction and the boundary
data $f$ and $h\in \mathfrak{D}'(\mathbb{T})$, the space of distributions in  $\mathbb{T}$.

Note that a solution to the biharmonic equation ${\Delta}^2 \Phi=0$ is called a {\it biharmonic function}.
See Almansi \cite{Al99}, Vekua \cite{Vekua} and \cite{AA1, AAK, AAK1, CPW0, CPW4} for properties of biharmonic functions.

The present article is motivated from the related studies in \cite{AKM2008, kalaj2008, KP, OL, P}.
In \cite{OL}, Olofsson considered the representation formula of the solutions to the following homogeneous
biharmonic Dirichlet problem (briefly, HBDP in the following):
\be\label{eq-11}
\begin{cases} \Delta^2 u=0 &
\mbox{ in  $\ID$},\\
\ds u=f & \mbox{ on $\mathbb{T}$},\\
\ds \partial_n u=h & \mbox{ on $\mathbb{T}$}.
\end{cases}
\ee
In particular, Olofsson proved the following (see \cite[Theorems 3.1 and 3.2]{OL}).

\begin{Thm}\label{pj-1} $(i)$ Suppose $u$ satisfies the growth conditions:
\beqq\label{eq-121}
|u(z)|\leq C(1-|z|)^{-N}\;\;\mbox{and}\;\;|\Delta u(z)|\leq C(1-|z|)^{-N},
\eeqq
where $C$ and $N$ are positive constants.
If $u$ satisfies HBDP \eqref{eq-11}, then $u$ admits the representation
\be\label{eq-12}
u(z)=F_0[f](z) +H_0[h](z),
\ee
where
\beqq
F_0[f](z)=\frac{1}{2\pi}\int^{2\pi}_{0}F_{0}(ze^{-i\theta})f(e^{i\theta})\,d\theta
~\mbox{ and }~
H_0[h](z)=\frac{1}{2\pi}\int^{2\pi}_{0}H_{0}(ze^{-i\theta})h(e^{i\theta})\,d\theta.
\eeqq
Here the kernels $H_0(z)$ and $F_0(z)$ are given by
\be\label{eq-1200}
H_0(z)=\frac{1}{2}\frac{(1-|z|^2)^2}{|1-z|^{2}} ~\mbox{ and }~F_0(z)=H_0(z)+\frac{1}{2}\frac{(1-|z|^2)^3}{|1-z|^{4}}.
\ee
\medskip
\noindent $(ii)$
If $u$ is defined by \eqref{eq-12}, then $u$ satisfies HBDP \eqref{eq-11}.
\end{Thm}

In fact, the function $F_0$ is a certain biharmonic Poisson kernel introduced by Abkar
and Hedenmalm \cite{AH}. Moreover, in \cite{BVZ, Be} the authors solved a certain Dirichlet boundary value problem for the polyharmonic equation in $\ID$.
In \cite{P}, Pavlovi\'{c} proved that the quasiconformality of harmonic homeomorphisms between $\ID$ can be characterized in terms
of their bi-Lipschitz continuity (\cite[Theorem 1.2]{P}). In \cite{AKM2008}, Arsenovi\'c et. al  showed that the Lipschitz continuity of
$\phi: \mathbb{S}^{n-1}\rightarrow\mathbb{R}^{n}$ implies the Lipschitz continuity of its harmonic extension
$P[\phi]:\mathbb{B}^{n}\rightarrow\mathbb{R}^{n}$ provided that $P[\phi]$ is a $K$-quasiregular mapping (\cite[Theorem 1]{AKM2008}), where
$\IB^n$ (resp. $\mathbb{S}^{n-1}$) denotes the unit ball (resp. the boundary of $\IB^n$) in $\IR^n$ and
$P$ stands for the usual Poisson kernel with respect to $\Delta$.
The assumption ``$P[\phi]$ being $K$-quasiregular" in \cite[Theorem 1]{AKM2008} is necessary as \cite[Example 1]{AKM2008} demonstrates.
Meanwhile, by assuming that $P[\phi]:\mathbb{B}^{n}\rightarrow\mathbb{B}^{n}$ is a $K$-quasiconformal harmonic mapping with $P[\phi](0)=0$ and
$\phi\in C^{1,\alpha}$, Kalaj \cite[Theorem 2.1]{kalaj2008} also proved the Lipschitz continuity of $P[\phi]$.

In particular, in \cite{KP}, Kalaj and Pavlovi\'{c} discussed the bi-Lipschitz continuity of quasiconformal
self-mappings in $\ID$ satisfying Poisson's equation $\Delta u=\psi$ (\cite[Theorem 1.2]{KP}).
See \cite{AH, AKM2008, clw, chen, K1, kalaj2008, KM0, KM1, KM2, KP0, MKM, LCW, P, MV} and the
references therein for detailed discussions on this topic.

In order to state our results, we need the representation formula of the biharmonic Green function in $\ID$.
The biharmonic Green function $G$ is the solution to the boundary value problem:
\beq\label{pj-3}
\begin{cases}
\Delta_z^2 \big(G(z, \zeta)\big)=\delta_{\zeta}(z) &  \mbox{for $z$ in $\ID$},\\
\ds G(z, \zeta)=0 & \mbox{for $z$ in $\mathbb{T}$},\\
\ds \frac{\partial G}{\partial n_z}(z, \zeta)=0 & \mbox{for $z$ in $\mathbb{T}$}
\end{cases}
\eeq
for each $\zeta\in\ID$, where $\delta_{\zeta}(z)$ denotes the Dirac distribution concentrated at the point $\zeta\in\ID$ and $\partial/\partial n_z$
stands for the inward normal derivative with respect to the variable $z\in\ID$. In $\ID$, the biharmonic Green function $G$ is given
by (cf. \cite{Al99, BVZ, Be})
\be\label{eq-150}
G(z, \zeta)=|z-\zeta|^2\log\left|\frac{1-\overline{\zeta}z}{z-\zeta}\right|^2-(1-|z|^2)(1-|\zeta|^2).
\ee

For convenience, we let
$$G[g](z)=\int_{\ID}G(z, \zeta)g(\zeta)\,d A(\zeta),
$$
where $dA(\zeta) =(1/\pi)\, dx\,dy$ denotes the normalized area measure in $\ID$.


Our first objective of this paper is to establish a representation formula and uniqueness of the solutions to IBDP \eqref{eq-16}, which is as follows.

\begin{thm}\label{thm-1.1}
$(1)$
Suppose that $u$ satisfies HBDP \eqref{eq-11}
and that
$w\in C^{\infty}(\overline{\ID})$ satisfies the following IBDP:
\be\label{eq-15}
\begin{cases}
\Delta^2 w=g & \mbox{ in $\ID$},\\
\ds w=0 & \mbox{ on $\mathbb{T}$},\\
\ds \partial_n w=0 & \mbox{ on $ \mathbb{T}$}.
\end{cases}
\ee
Then $\Phi=u+w$ is the only solution to IBDP \eqref{eq-16}, where
$$
u(z)=F_0[f](z)+H_0[h](z)\;\;\mbox{and}\;\; w(z)=G[g](z).
$$
\noindent $(2)$ If $g\in C(\overline{\ID})$,  then $\Phi$ solves IBDP \eqref{eq-16}, where
\be\label{pj-13}
\Phi(z)=F_0[f](z)+H_0[h](z)-G[g](z).
\ee
\end{thm}

The second objective is to discuss the bi-Lipschitz continuity of solutions to IBDP \eqref{eq-16}, which is formulated in the following form.

\begin{thm}\label{thm-1.2}
Suppose that $\Phi$ has the representation formula \eqref{pj-13}, $h\in C(\mathbb{T})$, $g\in C(\overline{\ID})$ and that $f$ satisfies the Lipschitz condition:
$$|f(e^{i\theta})-f(e^{i\varphi})|\leq L |e^{i\theta}-e^{i\varphi}|,
$$
where $L$ is a constant. Then for $z_1, z_2\in\ID$,
\be\label{eq-P-e1}
P\Big(\frac{Q}{P^2}-2\Big)|z_1-z_2| \leq |\Phi(z_1)-\Phi(z_2)|\leq P|z_1-z_2|,
\ee
where
\be\label{eq-P}
P=\frac{220}{3}L+4\|h\|_{\infty,\;\mathbb{T}}+\frac{23}{3}\|g\|_{\infty},
\ee
$$\|h\|_{\infty,\;\mathbb{T}}=\sup\{|h(z)|: z\in\mathbb{T}\},\;\; \|g\|_{\infty}=\sup\{|g(z)|: z\in\ID\},
$$
$$A=\left|\frac{1}{4\pi}\int^{2\pi}_{0}e^{-i\theta}\big(3f(e^{i\theta})+h(e^{i\theta})\big)\,d\theta
-\int_{\ID}\overline{\zeta}(\log|\zeta|^2+1-|\zeta|^2)g(\zeta)\,d A(\zeta)\right|^2,
$$
$$B=\left|\frac{1}{4\pi}\int^{2\pi}_{0}e^{i\theta}\big(3f(e^{i\theta})+h(e^{i\theta})\big)\,d\theta
-\int_{\ID}\zeta(\log|\zeta|^2+1-|\zeta|^2)g(\zeta)\,d A(\zeta)\right|^2
$$
and
 $Q=A-B$.
\end{thm}

It is worth pointing out that if $Q>2P^2$, then the condition \eqref{eq-P-e1} shows that $\Phi$ is bi-Lipschitz, otherwise,
$\Phi$ is Lipschitz. We would like to point out that this Lipschitz extension property is indeed interesting and does not hold
true for the classical Poisson kernel
$$P(z)=\frac{1-|z|^2}{|1-z|^{2}}.
$$

The rest of this article is organized as follows. In Section \ref{sec-2}, some necessary notations will be introduced and
several useful lemmas will be proved. In Section \ref{sec-3},  Theorem \ref{thm-1.1} will be
proved with the aid of Theorem \Ref{pj-1}. Section \ref{sec-4} will be devoted to the proof of Theorem \ref{thm-1.2}.

\section{Preliminaries}\label{sec-2}
In order to prove Theorems \ref{thm-1.1} and  \ref{thm-1.2} we need some preparation.


\subsection{Matrix norm}

Let
$$M=\left(
  \begin{array}{cc}
     a & b \\
     c & d \\
  \end{array}
\right)\in \mathbb{R}^{2\times2}.
$$
We will consider the matrix norm
$$\|M\|=\sup\{|Mz|:\; z\in \mathbb{C},\;|z|=1\}
$$
and the matrix function
$$l(M)=\inf\{|Mz|:\;z\in \mathbb{C},\;|z|=1\}.
$$
Let $D$ and $\Omega $ be plane domains in $\IC$. With $p(z) = u(z)+iv(z)$, $z = x+iy$, and $p:\,D\to \Omega$,  we can express the
Jacobian matrix $\nabla p$ of $p$ and Jacobian (determinant) $J(p)$ as
$$\nabla p =
\left(
  \begin{array}{cc}
    u_x & u_y \\
    v_x & v_y \\
  \end{array}
\right) ~\mbox{ and }~ J(p)=u_xv_y-v_xu_y.
$$
Obviously
\be\label{pj-11}
\|\nabla p(z)\|=\sup\{|\nabla p(z)\varsigma|:\; |\varsigma|=1\}=|p_z(z)|+|p_{\overline z}(z)|
\ee
and
$$
l(\nabla p(z))=\inf\{|\nabla p(z)\varsigma|:\; |\varsigma|=1\}=\big ||p_z(z)|-|p_{\overline z}(z)|\big |,
$$
where
$$p_{z}=\frac{\partial p}{\partial z}=\frac{1}{2}\left ( \frac{\partial p}{\partial x}-i\frac{\partial p}{\partial y}\right ) ~\mbox{ and }~
p_{\overline{z}}=\frac{\partial p}{\partial\overline{z}} =\frac{1}{2}\left ( \frac{\partial p}{\partial x}+i\frac{\partial p}{\partial y}\right ).
$$

We denote by $C_0^{\infty}(\ID)$ the space of functions which are infinitely differentiable and have compact support in $\ID$.
And $L_{loc}^1(\ID)$ denotes the space of locally integrable functions in $\ID$.

Functions $\psi\in L_{loc}^1(\ID)$ with distributions in $\ID$ having the action
$$\langle\psi, \varphi\rangle=\int_{\ID}\psi\varphi \,dA,\;\;\varphi\in C_0^{\infty}(\ID).
$$
Explicitly, we mean that the distribution $\psi_z$ have the action
$$\langle\psi_z, \varphi\rangle=-\int_{\ID}\psi\varphi_z \,dA,\;\;\varphi\in C_0^{\infty}(\ID),
$$
and similarly for the distribution $\psi_{\overline{z}}$ (cf. \cite{B}).

\subsection{Auxiliary results}

The following result is easy to derive: For $\beta>0$, we have
$$\frac{1}{2\pi}\int^{2\pi}_{0}\frac{d\theta}{|1-ze^{i\theta}|^{2\beta}}
=\sum^{\infty}_{n=0}\left(\frac{\Gamma(n+\beta)}{n!\Gamma(\beta)}\right)^2|z|^{2n},  \quad z\in\ID,
$$
where $\Gamma$ denotes the Gamma function.
Indeed, for $z\in\ID$, $|\zeta|=1$, and $\beta>0$, one has
$$\frac{1}{(1-z\zeta)^{\beta}}=\sum^{\infty}_{n=0}a_n z^n\zeta^n, \quad a_n=\frac{\Gamma(n+\beta)}{n!\Gamma(\beta)},
$$
and thus, by Parseval's theorem, we get
$$\frac{1}{2\pi}\int^{2\pi}_{0}\frac{d\theta}{|1-ze^{i\theta}|^{2\beta}}=
\frac{1}{2\pi}\int^{2\pi}_{0}\frac{d\theta}{|(1-ze^{i\theta})^{\beta}|^{2}}
=\sum^{\infty}_{n=0}|a_n|^2|z|^{2n}
$$
as required.

In particular, for  $\beta=1$,  $\beta=2$ and $\beta=3$, we obtain
\be\label{pj-14}
\frac{1}{2\pi}\int^{2\pi}_{0}\frac{d\theta}{|1-\overline{z}e^{-i\theta}|^{2}}=\sum^{\infty}_{n=0}|z|^{2n}=\frac{1}{1-|z|^2},
\ee
\be\label{pj-14'}
\frac{1}{2\pi}\int^{2\pi}_{0}\frac{d\theta}{|1-\overline{z}e^{-i\theta}|^{4}}=\sum^{\infty}_{n=0}(n+1)^2|z|^{2n}=\frac{1+|z|^2}{(1-|z|^2)^3}
\ee
and
\be\label{pj-14-ex1}
\frac{1}{2\pi}\int^{2\pi}_{0}\frac{1}{|1-\overline{z}e^{-i\theta}|^{6}}\,d\theta = \sum^{\infty}_{n=0}\frac{(n+1)^2(n+2)^2}{4}|z|^{2n}
\ee
Note that the above identities also hold when $\overline{z}$ is replaced by $z$. Thus, by the H\"{o}lder inequality,
\eqref{pj-14} and \eqref{pj-14'}, we easily have
\beq\label{pj-14a}
\frac{1}{2\pi}\int^{2\pi}_{0}\frac{d\theta}{|1-z e^{-i\theta}|^{3}}
&\leq&\left(\frac{1}{2\pi}\int^{2\pi}_{0}\frac{1}{|1-z e^{-i\theta}|^{2}}\,d\theta\right)^{\frac{1}{2}}
\left(\frac{1}{2\pi}\int^{2\pi}_{0}\frac{1}{|1-z e^{-i\theta}|^{4}}\,d\theta\right)^{\frac{1}{2}}\nonumber\\
&=&\frac{\sqrt{1+|z|^2}}{(1-|z|^2)^2}
\eeq
and this will be also used in the proof of Lemma \ref{lem4}.

\subsection{Useful lemmas}
\begin{lem}\label{lem-2}
Let $F_0$ be given by \eqref{eq-1200}. Then for $z\in\ID$, we have
$$\frac{1}{2\pi}\int^{2\pi}_{0}F_0(ze^{-i\theta})\,d\theta=1.
$$
\end{lem}
\bpf
By \eqref{pj-14} and \eqref{pj-14'}, we obtain that
\beqq
\frac{1}{2\pi}\int^{2\pi}_{0}F_0(ze^{-i\theta})\,d\theta
&=&\frac{1}{2\pi}\int^{2\pi}_{0}\frac{(1-|z|^2)^2 }{2|1-\overline{z}e^{i\theta}|^{2}}\,d\theta
+\frac{1}{2\pi}\int^{2\pi}_{0}\frac{(1-|z|^2)^3}{2|1-\overline{z}e^{i\theta}|^{4}}\,d\theta\\
&=&\frac{1-|z|^2}{2}+\frac{1+|z|^2}{2}=1
\eeqq
and the desired conclusion follows.
\epf

\begin{lem}\label{pj-100}
Let $H_0$ and $F_0$ be given by \eqref{eq-1200}, and $G$ be defined by \eqref{eq-150}.
\begin{enumerate}
\item[{\rm (a)}]  For any $\theta\in[0, 2\pi]$,
$$\frac{\partial H_0(ze^{-i\theta})}{\partial z}
=\frac{(1-|z|^2)[e^{-i\theta}(1-|z|^2)-2\overline{z}(1-ze^{-i\theta})]}{2(1-\overline{z}e^{i\theta})(1-ze^{-i\theta})^2}
$$
and
$$\frac{\partial H_0(ze^{-i\theta})}{\partial \overline{z}}=\overline{\left(\frac{\partial H_0(ze^{-i\theta})}{\partial z}\right)};
$$

\item[{\rm (b)}] For any $\theta\in[0, 2\pi]$,
\beqq\label{pj-12}
\frac{\partial F_0(ze^{-i\theta})}{\partial z}
&=&\frac{(1-|z|^2)[e^{-i\theta}(1-|z|^2)-2\overline{z}(1-ze^{-i\theta})]}{2(1-\overline{z}e^{i\theta})(1-ze^{-i\theta})^2}\\ \nonumber
&&+\frac{(1-|z|^2)^2[2e^{-i\theta}(1-|z|^2)-3\overline{z}(1-ze^{-i\theta})]}{2(1-\overline{z}e^{i\theta})^2(1-ze^{-i\theta})^3}
\eeqq
and
$$\frac{\partial F_0(ze^{-i\theta})}{\partial \overline{z}}=\overline{\left(\frac{\partial F_0(ze^{-i\theta})}{\partial z}\right)};
$$

\item[{\rm (c)}] For any fixed $\zeta\in \ID$,
$$G_z(z, \zeta)
=(\overline{z}-\overline{\zeta})\log\left|\frac{z-\zeta}{1-\overline{\zeta}z}\right|^2
+\frac{(\overline{z}-\overline{\zeta})(1-|\zeta|^2)}{1-\overline{\zeta}z}
-\overline{z}(1-|\zeta|^2)
$$
and
$$G_{\overline{z}}(z, \zeta)=\overline{G_z(z, \zeta)}.
$$
\end{enumerate}
\end{lem}

\begin{lem}\label{lem-4.1}
Suppose that $f, h \in C(\mathbb{T})$. Then
\beq\label{pj-401}
\frac{\partial F_0[f](z)}{\partial z}
=\frac{1}{2\pi}\int^{2\pi}_{0}\frac{\partial F_0(ze^{-i\theta})}{\partial z}f(e^{i\theta})\,d\theta,
\eeq
$$
\frac{\partial F_0[f](z)}{\partial \overline{z}}
=\frac{1}{2\pi}\int^{2\pi}_{0}\frac{\partial F_0(ze^{-i\theta})}{\partial \overline{z}}f(e^{i\theta})\,d\theta,
$$
$$
\frac{\partial H_0[h](z)}{\partial z}=\frac{1}{2\pi}\int^{2\pi}_{0}\frac{\partial H_0(ze^{-i\theta})}{\partial z}h(e^{i\theta})\,d\theta
$$
and
$$
\frac{\partial H_0[h](z)}{\partial \overline{z}}
=\frac{1}{2\pi}\int^{2\pi}_{0}\frac{\partial H_0(ze^{-i\theta})}{\partial \overline{z}}h(e^{i\theta})\,d\theta.
$$
\end{lem}
\bpf By Lemma \ref{pj-100}, we see that the functions
$$F_0(ze^{-i\theta})f(e^{i\theta}),\;\;\frac{\partial F_0(ze^{-i\theta})}{\partial z}f(e^{i\theta})\;\;\mbox{and}\;\;
\frac{\partial F_0(ze^{-i\theta})}{\partial \overline{z}}f(e^{i\theta})
$$
are all continuous on $\overline{\ID}_r\times[0, 2\pi]$, where $r\in[0, 1)$.

Let $z=\rho e^{i\varphi}\in \overline{\ID}_r$. Then we have
\be\label{pj-404a}
\frac{\partial F_0(ze^{-i\theta})}{\partial \rho}=\frac{1}{\rho} \left (\frac{\partial F_0(ze^{-i\theta})}{\partial z}z
+ \frac{\partial F_0(ze^{-i\theta})}{\partial \overline{z}}\overline{z}\right )
\ee
and
\be\label{pj-404b}
\frac{\partial F_0(ze^{-i\theta})}{\partial \varphi}=i\left (\frac{\partial F_0(ze^{-i\theta})}{\partial z}z
- \frac{\partial F_0(ze^{-i\theta})}{\partial \overline{z}}\overline{z}\right )
\ee
so that
\be\label{pj-404c}
\frac{\partial F_0(ze^{-i\theta})}{\partial z}=\frac{e^{-i\varphi}}{2}\left(\frac{\partial
F_0(ze^{-i\theta})}{\partial \rho}-\frac{i}{\rho}\frac{\partial F_0(ze^{-i\theta})}{\partial \varphi}\right).
\ee

It follows from \eqref{pj-404a} and \eqref{pj-404b} that both
$$\frac{\partial F_0(ze^{-i\theta})}{\partial \rho}f(e^{i\theta})\;\;\mbox{ and} \;\;\frac{\partial F_0(ze^{-i\theta})}{\partial \varphi}f(e^{i\theta})
$$
are continuous in $\overline{\ID}_r\times[0, 2\pi]$. Hence
\beqq
\int^{\rho}_{0}\int^{2\pi}_{0}\frac{\partial F_0(ze^{-i\theta})}{\partial \rho}f(e^{i\theta})\,d\theta\, d\rho
&=&\int^{2\pi}_{0}\int^{\rho}_{0}\frac{\partial F_0(ze^{-i\theta})}{\partial \rho}f(e^{i\theta})\,d\rho \,d\theta\\
&=&\int^{2\pi}_{0}\big(F_0(ze^{-i\theta})-F_0(0)\big)f(e^{i\theta})\,d\theta.
\eeqq
By differentiating with respect to $\rho$, we get
\beq\label{pj-405}
\int^{2\pi}_{0}\frac{\partial F_0(ze^{-i\theta})}{\partial \rho}f(e^{i\theta})\,d\theta
=\frac{\partial}{\partial \rho}\int^{2\pi}_{0}F_0(ze^{-i\theta})f(e^{i\theta})\,d\theta.
\eeq
On the other hand,
\beqq
\int^{\varphi}_{0}\int^{2\pi}_{0}\frac{\partial F_0(ze^{-i\theta})}{\partial \varphi}f(e^{i\theta})\,d\theta \,d\varphi
&=&\int^{2\pi}_{0}\int^{\varphi}_{0}\frac{\partial F_0(ze^{-i\theta})}{\partial \varphi}f(e^{i\theta})\,d\varphi \,d\theta\\
&=&\int^{2\pi}_{0}\big(F_0(ze^{-i\theta})-F_0(\rho e^{-i\theta})\big)f(e^{i\theta})\,d\theta.
\eeqq
By differentiating with respect to $\varphi$, we get
\beq\label{pj-406}
\int^{2\pi}_{0}\frac{\partial F_0(ze^{-i\theta})}{\partial \varphi}f(e^{i\theta})\,d\theta
=\frac{\partial}{\partial \varphi}\int^{2\pi}_{0}F_0(ze^{-i\theta})f(e^{i\theta})\,d\theta.
\eeq
It follows from  \eqref{pj-404c}, \eqref{pj-405}, and \eqref{pj-406} that \eqref{pj-401} holds.
The remaining results of the lemma follows similarly.
\epf

\begin{lem} \label{2.40}
Suppose that $F$ is defined on $\ID\times\ID$ and satisfies the following conditions:
\ben
\item[{\rm (a)}] $\frac{\partial F(z, \zeta)}{\partial z}$ and $\frac{\partial F(z, \zeta)}{\partial \overline{z}}$ exist;
\item[{\rm (b)}] For all $z\in\ID$, $\int_{\ID}\left|F(z, \zeta)\right|d A(\zeta)<\infty$;
\item[{\rm (c)}] For all $z\in\overline{\ID_r}$ $(r\in(0, 1))$
\beq\label{eq-2.40}
\int_{\ID}\left|\frac{\partial F(z, \zeta)}{\partial z}\right|d A(\zeta)<\infty\;
\mbox{and }\; \int_{\ID}\left|\frac{\partial F(z, \zeta)}{\partial \overline{z}}\right|d A(\zeta)<\infty.
\eeq
\een
Then for every $z\in\ID$, we have
\be\label{pj-20}
\frac{\partial }{ \partial z}\int_{\ID}F(z, \zeta)d A(\zeta)=\int_{\ID}\frac{\partial F(z, \zeta)}{\partial z}\,d A(\zeta)
\ee
and
\be\label{pj-21}
\frac{\partial }{ \partial \overline{z}}\int_{\ID}F(z, \zeta)d A(\zeta)=\int_{\ID}\frac{\partial F(z, \zeta)}{\partial \overline{z}}\,d A(\zeta).
\ee
\end{lem}
\bpf Obviously, we only need to prove the equality \eqref{pj-20} since the proof of \eqref{pj-21} is similar.
Let $z=\rho e^{i\varphi}\in\ID$. Then, \eqref{pj-404a} and \eqref{pj-404b} continue to hold if  $F_0(ze^{-i\theta})$
is replaced by $F(z, \zeta)$.
Thus, \eqref{eq-2.40} implies
$$\int^{\rho}_{0}\int_{\ID}\left|\frac{\partial F(z, \zeta)}{\partial \rho}\right|d A(\zeta)\,d\rho<\infty\;\;
\mbox{and }\;\;
\int^{\varphi}_{0}\int_{\ID}\left|\frac{\partial F(z, \zeta)}{\partial \varphi}\right|d A(\zeta)\,d\varphi<\infty.
$$
Now Fubini's theorem guarantees that
\beqq
\int^{\rho}_{0}\int_{\ID}\frac{\partial F(z, \zeta)}{\partial \rho}\,d A(\zeta)\,d \rho
=\int_{\ID}\int^{\rho}_{0}\frac{\partial F(z, \zeta)}{\partial \rho}\,d \rho \,d A(\zeta)
=\int_{\ID}(F(z, \zeta)-F(0, \zeta))\,d A(\zeta).
\eeqq
By differentiating with respect to $\rho$, we get
$$
\int_{\ID}\frac{\partial F(z, \zeta)}{\partial \rho}\,d A(\zeta)=\frac{\partial }{ \partial \rho}\int_{\ID}F(z, \zeta) \,d A(\zeta).
$$
Similarly, we can obtain that
$$
\int_{\ID}\frac{\partial F(z, \zeta)}{\partial \varphi}\,d A(\zeta)=\frac{\partial }{ \partial \varphi}\int_{\ID}F(z, \zeta)\,d A(\zeta).
$$
Since \eqref{pj-404c} continues to hold with $F(z, \zeta)$ in place $F_0(ze^{-i\theta})$,
\eqref{pj-20} obviously holds, and thus, the proof of the lemma is complete.
\epf

Recall that for $p\geq 1$ and $a>-1$, we have
\be\label{polylog1}
\int^{1}_{0}t^{a}\left (\log\frac{1}{t}\right )^{p-1} \,dt=\frac{\Gamma (p)}{(1+a)^p}
\ee
and the change of variable $t=r^2$ gives
\be\label{polylog2}
\int^{1}_{0}r^{2a+1}\left (\log\frac{1}{r^2}\right )^{p-1}\, dr=\frac{\Gamma (p)}{2(1+a)^p}.
\ee

\begin{lem} \label{lem-2.6}
For $g\in C(\overline{\ID})$, we have the following:
\ben
\item $\ds \int_{\ID}\left|G(z, \zeta)g(\zeta)\right|d A(\zeta)\leq \frac{3}{4}\|g\|_{\infty};$
\item
\ben
\item $\ds \frac{\partial G[g](z)}{\partial z}=\int_{\ID}\frac{\partial G(z, \zeta)}{\partial z}g(\zeta)\,d A(\zeta)$,
and
\item $\ds \left| \frac{\partial G[g](z)}{\partial z}\right|\leq \int_{\ID}\left|G_{z}(z, \zeta)g(\zeta)\right|d A(\zeta)
\leq \frac{23}{6}\|g\|_{\infty};$
\een
\item
\ben
\item
$\ds \frac{\partial G[g](z)}{\partial \overline{z}}=\int_{\ID}\frac{\partial G(z, \zeta)}{\partial \overline{z}}g(\zeta)\,d A(\zeta)$,
and
\item
$\ds \left| \frac{\partial G[g](z)}{\partial \overline{z}}\right|\leq\int_{\ID}\left|G_{\overline{z}}(z, \zeta)g(\zeta)\right|d A(\zeta)
\leq \frac{23}{6}\|g\|_{\infty}.$
\een\een
\end{lem}
\bpf
It follows from Lemma \ref{pj-100}(c) that for $z\in\ID$,
$$\int_{\ID} \left|G_z(z, \zeta)\right|d A(\zeta)\leq J_1+J_2+J_3,
$$
where
$$J_1=\int_{\ID}|z-\zeta|\log\left|\frac{1-\overline{\zeta}z}{z-\zeta}\right|^2\,d A(\zeta),
~J_2=\int_{\ID}\frac{(1-|\zeta|^2)|z-\zeta|}{|1-\overline{\zeta}z|}\,d A(\zeta)
$$
and
$$J_3=\int_{\ID}|z|(1-|\zeta|^2)\,d A(\zeta).
$$
Next, we estimate $J_1$, $J_2$ and $J_3$, respectively. In order to estimate $J_1$, we let
$$\zeta \mapsto \eta =\phi (\zeta)=\frac{z-\zeta}{1-\zeta\overline{z}}=re^{i\theta}
$$
so that $\phi =\phi^{-1}$,
$$\zeta=\frac{z-\eta}{1-\eta\overline{z}},~~z-\zeta=\frac{\eta(1-|z|^2)}{1-\eta\overline{z}}, ~~\phi ' (\zeta)=-\frac{1-|z|^2}{(1-\zeta\overline{z})^2},
$$
and thus,
$$d A(\zeta)= |({\phi^{-1}}) ' (\eta)|^2 d A(\eta)=\frac{(1-|z|^2)^2}{|1-\eta\overline{z}|^4}d A(\eta).
$$
Consequently, switching to polar coordinates yields
\beqq
J_1&=&\int_{\ID}\frac{|\eta|(1-|z|^2)^3}{|1-\eta\overline{z}|^5}\log\frac{1}{|\eta|^2}\,d A(\eta)
=\frac{(1-|z|^2)^3}{\pi}\int^{1}_{0}\int^{2\pi}_{0}\frac{r^2}{|1-\overline{z}re^{i\theta}|^5}\log\frac{1}{r^2}\,d\theta\, dr.
\eeqq
By the H\"{o}lder inequality,  \eqref{pj-14'} and \eqref{pj-14-ex1}, we get
\beqq
\frac{1}{2\pi}\int^{2\pi}_{0}\frac{d\theta}{|1-\overline{z}re^{i\theta}|^{5}}
&\leq&\left(\frac{1}{2\pi}\int^{2\pi}_{0}\frac{1}{|1-\overline{z}re^{i\theta}|^{4}}\,d\theta\right)^{\frac{1}{2}}
\left(\frac{1}{2\pi}\int^{2\pi}_{0}\frac{1}{|1-\overline{z}re^{i\theta}|^{6}}\,d\theta\right)^{\frac{1}{2}}\\
&=&\left(\sum^{\infty}_{n=0}(n+1)^2|z|^{2n}r^{2n}\right)^{\frac{1}{2}}\left(\sum^{\infty}_{n=0}\frac{(n+1)^2(n+2)^2}{4}|z|^{2n}r^{2n}\right)^{\frac{1}{2}}\\
&\leq&\sum^{\infty}_{n=0}\frac{(n+1)^2(n+2)^2}{4}|z|^{2n}r^{2n}
\eeqq
so that
\beqq
J_1&\leq&2(1-|z|^2)^3\sum^{\infty}_{n=0}\frac{(n+1)^2(n+2)^2}{4}|z|^{2n}\int^{1}_{0}r^{2n+2}\log\frac{1}{r^2} \,dr\\
&=&2(1-|z|^2)^3\sum^{\infty}_{n=0}\frac{(n+1)^2(n+2)^2}{2(2n+3)^2}|z|^{2n} \quad (\mbox{by \eqref{polylog2}})\\
&\leq&2(1-|z|^2)^3\sum^{\infty}_{n=0}\frac{(n+1)(n+2)}{8}|z|^{2n} =\frac{1}{2}
\eeqq
because  $4(n+1)(n+2)\leq(2n+3)^2$ and $\sum^{\infty}_{n=0}(n+1)(n+2)r^{n} =2/(1-r)^3$.
By the triangle inequality, we get
$$J_2\leq\int_{\ID}\frac{(1-|\zeta|^2)(1+|\zeta|)}{1-|\zeta|}\,d A(\zeta)
=\frac{1}{\pi}\int^{2\pi}_{0}\int^{1}_{0}(1+\rho)^2\rho \,d\rho \,d\varphi=\frac{17}{6},
$$
where $\zeta=\rho e^{i\varphi}$. Finally, we obtain
$$J_3\leq\frac{1}{\pi}\int^{2\pi}_{0}\int^{1}_{0}(1-\rho^2)\rho \,d\rho \,d\varphi=\frac{1}{2}
$$
The bounds on $J_1$, $J_2$ and $J_3$ give $\int_{\ID}\left|G_z(z, \zeta)\right|d A(\zeta)\leq 23/6$
and thus,
\beqq\label{eq-2.41}
\int_{\ID}\left|G_{z}(z, \zeta)g(\zeta)\right|d A(\zeta)\leq\frac{23}{6}\|g\|_{\infty} ~\mbox{ and }~
\int_{\ID}\left|G_{z}(z, \zeta)g(\zeta)\right|d A(\zeta)\leq \frac{23}{6}\|g\|_{\infty},
\eeqq
where the second inequality above is a consequence of Lemma \ref{pj-100}.
Moreover, by the similar reasoning as above, we deduce that
$$\int_{\ID}\left|G(z, \zeta)\right|d A(\zeta) \leq I+\int_{\ID}(1-|\zeta|^2)\, d A(\zeta)=I+\frac{1}{2},
$$
where
\beq
I& =& \int_{\ID}|z-\zeta|^2\log\left|\frac{1-\overline{\zeta}z}{z-\zeta}\right|^2\, d A(\zeta) \nonumber\\
&=&\int_{\ID}\frac{|\eta|^2(1-|z|^2)^4}{|1-\eta\overline{z}|^6}\log\frac{1}{|\eta|^2} \,d A(\eta) \nonumber\\
&=&\frac{(1-|z|^2)^4}{\pi}\int^{1}_{0}\int^{2\pi}_{0}\frac{r^3}{|1-\overline{z}re^{i\theta}|^6}\log\frac{1}{r^2}\,d\theta \,dr\nonumber \\
&=&2(1-|z|^2)^4\sum^{\infty}_{n=0}\frac{(n+1)^2(n+2)^2}{4}|z|^{2n}\int^{1}_{0}r^{2n+3}\log\frac{1}{r^2} \,dr \nonumber \\
&=&(1-|z|^2)^4\sum^{\infty}_{n=0}\frac{(n+1)^2}{4}|z|^{2n}  \quad (\mbox{by \eqref{polylog2}})\nonumber\\
&=&(1-|z|^2)^4\left (\frac{1+|z|^2}{4(1-|z|^2)^3}\right ) =\frac{1-|z|^4}{4}\nonumber
\eeq
and thus, we have
$$\int_{\ID}\left|G(z, \zeta)\right|d A(\zeta)\leq I+\frac{1}{2}\leq\frac{3}{4}.
$$
Hence the proof of the lemma is complete, since the rest of it follows from Lemma \ref{2.40}.
\epf

\section{Solutions to IBDP \eqref{eq-16}} \label{sec-3}
First, let us recall a useful result from \cite{OL}.

\begin{Thm}\label{pj-10} $($\cite[Theorem 2.1]{OL}$)$
Suppose $u$ satisfies the conditions:
$$\Delta^2 u=0\;\;\mbox{and}\;\;\lim_{r\to 1}\frac{u(re^{i\theta})}{1-r}=0\;\;\mbox{in}\;\,\mathfrak{D}'(\mathbb{T}),
$$
where $re^{i\theta}\in \mathbb{D}$. Then $u=0$ in $\mathbb{D}$.
\end{Thm}

Our next result concerns the representation formula and the uniqueness of the solutions to the  IBDP \eqref{eq-15}.

\begin{lem}\label{lem3.1} The function $G(z, \zeta)$ given by \eqref{eq-150} satisfies the following:
\ben
\item[{\rm (a)}]
$\int_{\ID}\left|G(z, \zeta)\right|d A(z)<\infty$ and $\int_{\ID} \left|G_z(z, \zeta)\right|d A(z)<\infty;$
\item[{\rm (b)}]
For fixed $\zeta\in\ID$,
$$
G_{z\overline{z}}(z, \zeta)=H_2(z, \zeta)=\log\left|\frac{1-\overline{\zeta}z}{z-\zeta}\right|^2
-\frac{(1-|\zeta|^2)(1-|z|^2|\zeta|^2)}{|1-z\overline{\zeta}|^2}
$$
in the sense of distributions in $\ID;$
\item[{\rm (c)}]
$\int_{\ID}\left|G_{z\overline{z}}(z, \zeta)\right|d A(\zeta)<\infty$ and
$\int_{\ID}\left|G_{z\overline{z}}(z, \zeta)\right|d A(z)<\infty$.
\een
\end{lem}
\bpf
By \eqref{eq-150}, we find that
$$
\int_{\ID}\left|G(z, \zeta)\right|d A(z) \leq J_4+\int_{\ID}(1-|z|^2)\,d A(z)\leq J_4+\frac{1}{2},
$$
where
$$J_4=\int_{\ID}|z-\zeta|^2\log\left|\frac{1-\overline{\zeta}z}{z-\zeta}\right|^2\,d A(z) .
$$
In order to estimate $J_4$, we let
$$z \mapsto \tau =\psi (z)=\frac{\zeta -z}{1-\overline{\zeta}z}
$$
so that $\psi =\psi^{-1}$,
%
%
$$z=\frac{\zeta-\tau}{1-\overline{\zeta}\tau}\;\;\mbox{and}\;\;d A(z)=\frac{(1-|\zeta|^2)^2}{|1-\tau\overline{\zeta}|^4}d A(\tau),
$$
and consequently, as in the proof of Lemma \ref{lem-2.6}, we have
$$J_4=\int_{\ID}\frac{|\tau|^2(1-|\zeta|^2)^4}{|1-\tau\overline{\zeta}|^6}\log\frac{1}{|\tau|^2}\, d A(\tau)
=\frac{1-|\zeta|^4}{4}.
$$
Hence
$$\int_{\ID}\left|G(z, \zeta)\right|d A(z)\leq J_4+\frac{1}{2}\leq \frac{3}{4}.
$$
Moreover, by the similar reasoning as above, we find that
\beqq
\int_{\ID} \left|G_z(z, \zeta)\right|d A(z)
&\leq &\int_{\ID}|z-\zeta|\log\left|\frac{1-\overline{\zeta}z}{z-\zeta}\right|^2\, d A(z)+2(1-|\zeta|^2)\int_{\ID} \,d A(z)\\
&\leq&\frac{1}{2} +2(1-|\zeta|^2)\leq \frac{5}{2}.
\eeqq
Thus the assertion (a) in the lemma is true. Therefore,
$G_z(z, \zeta)\in L^1(\ID)$ so that its derivative has the action
$$\langle G_{z\overline{z}}(z, \zeta), \varphi(z)\rangle=-\int_{\ID} G_z(z, \zeta)\varphi_{\overline{z}}(z)\, d A(z), \;\;\varphi\in C_0^{\infty}(\ID).
$$
By Lebesgue's dominated convergence theorem we get
$$\int_{\ID} G_z(z, \zeta)\varphi_{\overline{z}}(z)\, d A(z)
=\lim_{\varepsilon\to0}\int_{\ID\setminus \mathbb{D}(\zeta, \varepsilon)} G_z(z, \zeta)\varphi_{\overline{z}}(z) \, d A(z).
$$
For any small $\varepsilon>0$, let $D_{\varepsilon}=D(\zeta, \varepsilon)$. Partial integration gives
$$\int_{\ID\setminus D_{\varepsilon}} G_z(z, \zeta)\varphi_{\overline{z}}(z) \,d A(z)=\int_{\partial D_{\varepsilon}} G_z(z, \zeta)\varphi(z)v(z)\, d s(z)
-\int_{\ID\setminus D_{\varepsilon}}G_{z\overline{z}}(z, \zeta)\varphi(z)\,d A(z),
$$
where $v$ is the unit outward normal of $\ID\setminus D_{\varepsilon}$, that is, the inward unit normal of $D_{\varepsilon}$ and $ds$ denotes the normalized
arc length measure. It follows from Lemma \ref{pj-100} that
$$\left|G_z(z, \zeta)\right|\leq (1+|\zeta|)\log\left|\frac{1-\overline{\zeta}z}{z-\zeta}\right|^2+2(1-|\zeta|^2)
\leq 2\left(1-\log\left|\frac{z-\zeta}{1-\overline{\zeta}z}\right|^2\right).
$$
Then
\beqq
\left|\int_{\partial D_{\varepsilon}} G_z(z, \zeta)\varphi(z)v(z)\, d s(z)\right|
&\leq& C_{\varphi}\int_{\partial D_{\varepsilon}} \left|G_z(z, \zeta)\right|d s(z)\\
&\leq& 2C_{\varphi}(1-\log {\varepsilon}^2)\varepsilon\rightarrow0
\eeqq
as $\varepsilon\rightarrow0$, where $C_{\varphi}$ is a constant depending only on $\sup\varphi$.
A straightforward computation shows that for $z\neq \zeta$ we have
$$G_{z\overline{z}}(z, \zeta)=H_2(z,\zeta)=\log\left|\frac{1-\overline{\zeta}z}{z-\zeta}\right|^2
-\frac{(1-|\zeta|^2)(1-|z|^2|\zeta|^2)}{|1-z\overline{\zeta}|^2}
$$
so that
\be\label{eq-sam1}
\left|G_{z\overline{z}}(z, \zeta)\right|\leq \log\left|\frac{1-\overline{\zeta}z}{z-\zeta}\right|^2+2(1+|\zeta|)+1-|\zeta|^2
\leq 5\left(1-\log\left|\frac{z-\zeta}{1-\overline{\zeta}z}\right|^2\right).
\ee
Therefore
\beqq
\langle G_{z\overline{z}}(z, \zeta), \varphi(z)\rangle
&=&\lim_{\varepsilon\to0}\int_{\ID\setminus D_{\varepsilon}}G_{z\overline{z}}(z, \zeta)\varphi(z)\, d A(z)\\
&=&\int_{\ID}H_2(z,\zeta)\varphi(z)\, d A(z)
=\langle H_2(z,\zeta), \varphi(z)\rangle.
\eeqq
The assertion (b) in the lemma is hold. To prove the assertion (c),  it follows from the assertion (b) and the first inequality in \eqref{eq-sam1} that
\beqq
\int_{\ID}\left|G_{z\overline{z}}(z, \zeta)\right|d A(\zeta)
&\leq&\int_{\ID}\log\left|\frac{1-\overline{\zeta}z}{z-\zeta}\right|^2\,d A(\zeta)+C.
\eeqq

We make a convention that in the course of the proof, the value of constants may change from one occurrence to the next,
but we always use the same letter $C$ to denote them.

It follows directly from Green representation formula that
$$\int_{\ID}\log\left|\frac{1-\overline{\zeta}z}{z-\zeta}\right|^2\,d A(\zeta)=1-|z|^2\leq1,
$$
which implies that
$\int_{\ID}\left|G_{z\overline{z}}(z, \zeta)\right|d A(\zeta)<\infty.
$
By the similar reasoning as above, one obtains that
$\int_{\ID}\left|G_{z\overline{z}}(z, \zeta)\right|d A(z)<\infty.$
The proof of the lemma is complete.
\epf

\begin{lem}\label{lem3.2}
The function $G(z, \zeta)$ satisfies:
\ben
\item[{\rm (a)}]
For fixed $\zeta\in\ID$,
$$G_{z\overline{z}z}(z, \zeta) =H_3(z, \zeta)
=-\frac{1-|\zeta|^2}{(z-\zeta)(1-\overline{\zeta}z)}-\frac{\overline{\zeta}(1-|\zeta|^2)}{(1-\overline{\zeta}z)^2}
$$
in the sense of distributions in $\ID;$
\item[{\rm (b)}]
$\int_{\ID}\left|G_{z\overline{z}z}(z, \zeta)\right|d A(\zeta)<\infty$ and for fixed
$\zeta\in\ID$, $\int_{\ID}\left|G_{z\overline{z}z}(z, \zeta)\right|d A(z)<\infty$.
\een
\end{lem}
\bpf
It follows from Lemma \ref{lem3.1}(c) that for fixed $\zeta\in\ID$, $G_{z\overline{z}}(z, \zeta)\in L^1(\ID)$. Then
$$\langle G_{z\overline{z}z}(z, \zeta), \varphi(z)\rangle
=-\int_{\ID} G_{z\overline{z}}(z, \zeta)\varphi_{z}(z)\,d A(z), \;\;\varphi\in C_0^{\infty}(\ID).
$$
By Lebesgue's dominated convergence theorem we get
$$\int_{\ID} G_{z\overline{z}}(z, \zeta)\varphi_{z}(z)\, d A(z)
=\lim_{\varepsilon\to0}\int_{\ID\setminus D_{\varepsilon}} G_{z\overline{z}}(z, \zeta)\varphi_{z}(z)\,d A(z).
$$
Partial integration gives
\beqq
\int_{\ID\setminus D_{\varepsilon}} G_{z\overline{z}}(z, \zeta)\varphi_{z}(z)\,d A(z)
& = &-\int_{\partial D_{\varepsilon}} G_{z\overline{z}}(z, \zeta)\varphi(z)v(z)\, d s(z)\\
&& \hspace{2cm} -\int_{\ID\setminus D_{\varepsilon}}G_{z\overline{z}z}(z, \zeta)\varphi(z)\,d A(z).
\eeqq
By Lemma \ref{lem3.1} and the second inequality in \eqref{eq-sam1}, we get
\beqq
\left|\int_{\partial D_{\varepsilon}} G_{z\overline{z}}(z, \zeta)\varphi(z)v(z)\, d s(z)\right|
&\leq& C_{\varphi}\int_{\partial D_{\varepsilon}} \left|G_{z\overline{z}}(z, \zeta)\right|d s(z)\\
&\leq& 5C_{\varphi}(1-\log {\varepsilon}^2)\varepsilon\rightarrow0
\eeqq
as $\varepsilon\rightarrow0$. Moreover, for $z\neq \zeta$ we have
$$G_{z\overline{z}z}(z, \zeta)=H_3(z, \zeta)
=-\frac{1-|\zeta|^2}{(z-\zeta)(1-\overline{\zeta}z)}-\frac{\overline{\zeta}(1-|\zeta|^2)}{(1-\overline{\zeta}z)^2}
$$
and thus,
\beqq
\langle G_{z\overline{z}z}(z, \zeta), \varphi(z)\rangle
&=&\lim_{\varepsilon\to0}\int_{\ID\setminus D_{\varepsilon}}G_{z\overline{z}z}(z, \zeta)\varphi(z)\, d A(z)\\
&=&\int_{\ID}H_3(z,\zeta)\varphi(z)\, d A(z)\\
&=&\langle H_3(z,\zeta), \varphi(z)\rangle.\eeqq
The assertion (a) in the lemma holds.

By Lemma \ref{lem3.2}(a),
$$\int_{\ID}\left|G_{z\overline{z}z}(z, \zeta)\right|\, d A(\zeta)
\leq\int_{\ID}\frac{1-|\zeta|^2}{|z-\zeta|\cdot|1-\overline{\zeta}z|}\, d A(\zeta)
+\int_{\ID}\frac{1-|\zeta|^2}{|1-\overline{\zeta}z|^2}\,d A(\zeta).
$$
Moreover, as before, the transformation
$$\zeta \mapsto \eta =\phi (\zeta)=\frac{z-\zeta}{1-\zeta\overline{z}}=re^{i\theta}
$$
gives after some computation that
\beqq
\int_{\ID}\frac{1-|\zeta|^2}{|z-\zeta|\cdot|1-\overline{\zeta}z|}\, d A(\zeta)
&=&\int_{\ID}\frac{(1-|z|^2)(1-|\eta|^2)}{|\eta|\cdot|1-\eta\overline{z}|^4}\, d A(\eta)\\
&=&\frac{(1-|z|^2)}{\pi}\int^{1}_{0}\int^{2\pi}_{0}\frac{1-r^2}{|1-\overline{z}re^{i\theta}|^4}\,d\theta\, dr\\
&=&2(1-|z|^2)\sum^{\infty}_{n=0}(n+1)^2|z|^{2n}\int^{1}_{0}r^{2n}(1-r^2)\, dr ~~(\mbox{by \eqref{pj-14'}})\\
&=&4(1-|z|^2)\sum^{\infty}_{n=0}\frac{(n+1)^2}{(2n+1)(2n+3)} |z|^{2n}  \\
&\leq&\frac{4(1-|z|^2)}{3}\sum^{\infty}_{n=0} |z|^{2n} = \frac{4}{3} 
\eeqq
and similarly, it is easy to see that
$$\int_{\ID}\frac{1-|\zeta|^2}{|1-\overline{\zeta}z|^2}\, d A(\zeta)
=\int_{\ID}\frac{(1-|z|^2)(1-|\eta|^2)}{|1-\eta\overline{z}|^4}\, d A(\eta)\leq1,
$$
Thus, we conclude that
\beq\label{eq-33}\nonumber
\int_{\ID}\left|G_{z\overline{z}z}(z, \zeta)\right|d A(\zeta)<\infty.
\eeq
On the other hand, since $|G_{z\overline{z}z}|\leq 2(1+|\zeta|)/(1-|\zeta|)$, we have
$$\int_{\ID}\left|G_{z\overline{z}z}(z, \zeta)\right|d A(z)\leq 2\frac{1+|\zeta|}{1-|\zeta|}<\infty
$$
for each fixed $\zeta\in\ID$. The proof of the lemma is complete.
\epf

\begin{lem}\label{pj-2}
\ben
\item[{\rm (a)}] If $w\in C^{\infty}(\overline{\ID})$ is a solution to \eqref{eq-15}, then
\beq\label{eq-31}
w(z)=G[g](z)=\int_{\ID}G(z, \zeta)g(\zeta)\,d A(\zeta),
\eeq
where $G(z, \zeta)$ is defined by \eqref{eq-150}.

\item[{\rm (b)}]
If $g\in C(\overline{\ID})$ and $w$ is defined by \eqref{eq-31}, then $w$ is the only solution to \eqref{eq-15}.
\een
\end{lem}
\bpf
Part (a) is a consequence of applying Green's formula twice  (cf. \cite{AA} or \cite{HJS})

To prove the second part, we assume that $w$ admits the
representation formula \eqref{eq-31}. First, we prove that $w$ satisfies the first equation in IBDP \eqref{eq-15}.
\bcl\label{cl-3}
$\Delta_z^2\Big(\int_{\ID}G(z, \zeta)g(\zeta)\,d A(\zeta)\Big) =g(z).$
\ecl
Obviously, it follows from \eqref{pj-3} that
$$\int_{\ID}\Delta_z^2 \big(G(z, \zeta)\big)g(\zeta)\,d A(\zeta)=g(z)
$$
and thus, it suffices to show that
$$\Delta_z^2\Big(\int_{\ID}G(z, \zeta)g(\zeta)\,d A(\zeta)\Big)=\int_{\ID}\Delta_z^2 \big(G(z, \zeta)\big)g(\zeta)\, d A(\zeta)
$$
in the sense of distributions in $\ID$.

To prove this equality, we divide the proof into four steps.
\begin{step}
By Lemma \ref{lem-2.6}, we obtain that
$$\frac{\partial}{\partial z}\Big(\int_{\ID}G(z, \zeta)g(\zeta)\,d A(\zeta)\Big)
=\int_{\ID}G_z(z, \zeta)g(\zeta)\, d A(\zeta).
$$
\end{step}
\begin{step}\label{step-2}
We prove that $\frac{\partial^2}{\partial z\partial \overline{z}}\Big(\int_{\ID}G(z, \zeta)g(\zeta)\, d A(\zeta)\Big)
=\int_{\ID}G_{z\overline{z}}(z, \zeta)g(\zeta)\, d A(\zeta)$ in the sense of distributions in $\ID$.
\end{step}
To prove the second step, we recall from Lemma \ref{lem-2.6} that
\beqq
\left\langle\frac{\partial^2}{\partial z\partial \overline{z}}\Big(\int_{\ID}G(z, \zeta)g(\zeta)\, d A(\zeta)\Big), \varphi(z)\right\rangle
&=&\left\langle\int_{\ID}G(z, \zeta)g(\zeta)\, d A(\zeta), \varphi_{z\overline{z}}(z)\right\rangle\\
&=&\int_{\ID}\left(\int_{\ID}G(z, \zeta)g(\zeta)\, d A(\zeta)\varphi_{z\overline{z}}(z)\right)d A(z),
\eeqq
where $\varphi\in C_0^{\infty}(\ID)$.
By Lemmas \ref{lem-2.6} and \ref{lem3.1}, we know that
\beqq
\int_{\ID}\left(\int_{\ID}G(z, \zeta)g(\zeta)\, d A(\zeta)\varphi_{z\overline{z}}(z)\right)d A(z)
&=&\int_{\ID}\left(\int_{\ID}G(z, \zeta)g(\zeta)\varphi_{z\overline{z}}(z)\, d A(z)\right)d A(\zeta)\\
&=&\int_{\ID}\left\langle G(z, \zeta)g(\zeta), \varphi_{z\overline{z}}(z)\right\rangle \,d A(\zeta)\\
&=&\int_{\ID}\left\langle G_{z\overline{z}}(z, \zeta)g(\zeta), \varphi(z)\right\rangle \,d A(\zeta)\\
&=&\int_{\ID}\left(\int_{\ID}G_{z\overline{z}}(z, \zeta)g(\zeta)\varphi(z)\,d A(z)\right)d A(\zeta)\\
&=&\int_{\ID}\left(\int_{\ID}G_{z\overline{z}}(z, \zeta)g(\zeta)\, d A(\zeta)\varphi(z)\right)d A(z)\\
&=&\left\langle\int_{\ID}G_{z\overline{z}}(z, \zeta)g(\zeta)\,d A(\zeta), \varphi(z)\right\rangle.
\eeqq
Hence
$$\left\langle\frac{\partial^2}{\partial z\partial \overline{z}}\Big(\int_{\ID}G(z, \zeta)g(\zeta)\,d A(\zeta)\Big), \varphi(z)\right\rangle
=\left\langle\int_{\ID}G_{z\overline{z}}(z, \zeta)g(\zeta)\, d A(\zeta), \varphi(z)\right\rangle
$$
as required. Hence, the proof of Step \ref{step-2} is finished.

\begin{step}\label{step-3}
We prove that $\frac{\partial^3}{\partial z\partial \overline{z}\partial z}\Big(\int_{\ID}G(z, \zeta)g(\zeta)\, d A(\zeta)\Big)
=\int_{\ID}G_{z\overline{z}z}(z, \zeta)g(\zeta)\, d A(\zeta)$ in the sense of distributions in $\ID$.
\end{step}

By the similar reasoning as in the proof of Step \ref{step-2}, it follows from
Lemmas \ref{lem-2.6}, \ref{lem3.1} and \ref{lem3.2} and Step \ref{step-2}, we see that
$$\left\langle\frac{\partial^3}{\partial z\partial \overline{z}z}\Big(\int_{\ID}G(z, \zeta)g(\zeta)\, d A(\zeta)\Big), \varphi(z)\right\rangle
=\left\langle\int_{\ID}G_{z\overline{z}z}(z, \zeta)g(\zeta)\,d A(\zeta), \varphi(z)\right\rangle,
$$
where $\varphi\in C_0^{\infty}(\ID)$. Hence, the proof of Step \ref{step-3} is finished.

\begin{step}\label{step-4}
We prove that $\frac{\partial^4}{\partial z\partial \overline{z}\partial z\partial \overline{z}}\Big(\int_{\ID}G(z, \zeta)g(\zeta)\,
d A(\zeta)\Big)
=\int_{\ID}G_{z\overline{z}z\overline{z}}(z, \zeta)g(\zeta)\, d A(\zeta)$ in the sense of distributions in $\ID$.
\end{step}
It follows from $\Delta_z^2 \big(G(z, \zeta)\big)=\delta_{\zeta}(z)$ that
$$\left|\int_{\ID}G_{z\overline{z}z\overline{z}}(z, \zeta)g(\zeta)\, d A(\zeta)\right|\leq\|g\|_{\infty}.
$$
By the similar reasoning as in the proof of Step \ref{step-2}, it follows from Lemmas \ref{lem-2.6},
\ref{lem3.1} and \ref{lem3.2}, Steps \ref{step-2} and \ref{step-3}, that
$$\left\langle\frac{\partial^4}{\partial z\partial \overline{z}z\overline{z}}\Big(\int_{\ID}G(z, \zeta)g(\zeta)\,d A(\zeta)\Big), \varphi(z)\right\rangle
=\left\langle\int_{\ID}G_{z\overline{z}z\overline{z}}(z, \zeta)g(\zeta)\,d A(\zeta), \varphi(z)\right\rangle,
$$
where $\varphi\in C_0^{\infty}(\ID)$. This completes the proof of Step \ref{step-4} and this confirms
Claim \ref{cl-3}.


Next, we check that $w|_{\mathbb{T}}=0$, i.e., the middle equation in IBDP \eqref{eq-15} holds.

\bcl\label{cl-1}\label{cl-3.2}
$\lim_{z\rightarrow e^{it}}w(z)=0.$
\ecl

Again, it follows from \eqref{pj-3} that
$$\int_{\ID}\lim_{z\rightarrow e^{it}}G(z, \zeta)g(\zeta)\,d A(\zeta)=0
$$
and thus, to prove the claim, it suffices to show that
$$\lim_{z\rightarrow e^{it}}\int_{\ID}G(z, \zeta)g(\zeta)\,d A(\zeta)=\int_{\ID}\lim_{z\rightarrow e^{it}}G(z, \zeta)g(\zeta)\,d A(\zeta).
$$
To show this, we use the Vitali theorem (cf. \cite[Theorem 26]{Ha} or \cite[p.~4050]{KP}) which asserts that
if $X$ is a measurable space with finite measure $\mu$ and that $p_n: X\to \IC$ is a sequence of functions such that
$$\lim_{n\to\infty}p_n(x)=p(x)\;a.e.\;\;\mbox{and}\;\;\sup_n\int_X|p_n|^q\, d\mu<\infty\;\;\mbox{for some}\;\; q>1,
$$
then
$$\lim_{n\to \infty}\int_X p_n\, d\mu=\int_X p\,d\mu.
$$

By the Vitali theorem and \eqref{pj-3}, we only need to demonstrate that
$$\sup_{z\in\ID}\int_{\ID}|G(z, \zeta)g(\zeta)|^2 \,d A(\zeta)<\infty.
$$

Applying the Minkowski inequality, we conclude that
$$\left(\int_{\ID}|G(z, \zeta)|^2\, d A(\zeta)\right)^{\frac{1}{2}} \leq I_1 +I_2,
$$
where
$$I_1=\left(\int_{\ID}|z-\zeta|^4\log^2\left|\frac{1-\overline{\zeta}z}{z-\zeta}\right|^2\, d A(\zeta)\right)^{\frac{1}{2}}
$$
and
$$I_2=\left(\int_{\ID}(1-|z|^2)^2(1-|\zeta|^2)^2\,d A(\zeta)\right)^{\frac{1}{2}}\leq C.
$$
Again, the transformation
$$\zeta \mapsto \eta =\phi (\zeta)=\frac{z-\zeta}{1-\zeta\overline{z}}=re^{i\theta}
$$
gives after some computation that
\beqq
I_1^2
&=&\int_{\ID}\frac{|\eta|^4(1-|z|^2)^6}{|1-\eta\overline{z}|^8}\log^2\frac{1}{|\eta|^2}\, d A(\eta)\\
&=&2(1-|z|^2)^6\sum^{\infty}_{n=0}\frac{(n+1)^2(n+2)^2(n+3)^2}{36}|z|^{2n}\int^{1}_{0}r^{2n+5}\log^2\frac{1}{r^2}\, dr\\
&=&(1-|z|^2)^6\sum^{\infty}_{n=0}\frac{(n+1)^2(n+2)^2}{18(n+3)}|z|^{2n} \quad (\mbox{by \eqref{polylog2}})\\
&\leq&(1-|z|^2)^6\sum^{\infty}_{n=0}\frac{(n+1)^2(n+2)}{18}|z|^{2n} \\ 
&=&\frac{1}{18}(1-|z|^2)^2(2+4|z|^2)\leq\frac{1}{3}.
\eeqq
Then
$$\sup_{z\in\ID}\int_{\ID}|G(z, \zeta)g(\zeta)|^2\, d A(\zeta)<\infty,
$$
since $g\in C(\overline{\ID})$, which is what we want.\medskip

To finish the proof, we have to show that $w$ satisfies the third equation in \eqref{eq-15}.
\bcl\label{cl-2}
$\lim_{z\rightarrow e^{it}}\partial_{n(z)}w(z)=0$.
\ecl

Once again, it following from \eqref{pj-3} that
$$\int_{\ID}\lim_{z\rightarrow e^{it}}\partial_{n(z)}G(z, \zeta)g(\zeta)\, d A(\zeta)=0.
$$

To prove this claim, the similar reasoning stated as in the beginning of the proof of Claim \ref{cl-3.2} shows that it
suffices to check the boundedness of the integral
\be\label{eq-sam2}
\int_{\ID}\left|\partial_{n(z)}G(z, \zeta)g(\zeta)\right|^2\, d A(\zeta).
\ee

By letting $z=re^{i\theta}$, we write
$\partial_{n(z)}G(z, \zeta)=G_z(z, \zeta)e^{i\theta}+G_{\overline{z}}(z, \zeta)e^{-i\theta}.
$
Since $G_{\overline{z}}(z, \zeta)=\overline{G_z(z, \zeta)}$,
to prove the boundedness of \eqref{eq-sam2}, we only need to show that
$$\sup_{z\in\ID}\int_{\ID}|G_z(z, \zeta)g(\zeta)|^2\, d A(\zeta)<\infty.
$$
By the expression for $G_z$ from Lemma \ref{pj-100}(c) and the Minkowski inequality, we have
$$\left(\int_{\ID}|G_z(z, \zeta)|^2\, d A(\zeta)\right)^{\frac{1}{2}}
\leq J +C, \quad
J^2=\int_{\ID}|z-\zeta|^2\log^2\left|\frac{1-\overline{\zeta}z}{z-\zeta}\right|^2\, d A(\zeta).
$$
Now, as before, by letting $\eta=\frac{z-\zeta}{1-\zeta\overline{z}}=r e^{i\varphi}$, we find that
\beqq
J^2
&=&\int_{\ID}\frac{|\eta|^2(1-|z|^2)^4}{|1-\eta\overline{z}|^6}\log^2\frac{1}{|\eta|^2}\, d A(\eta)\\
&=&2(1-|z|^2)^4\sum^{\infty}_{n=0}\frac{(n+1)^2(n+2)^2}{4}|z|^{2n}\int^{1}_{0}r^{2n+3}\log^2\frac{1}{r^2}\, dr\\
&=&\frac{(1-|z|^2)^4}{2}\sum^{\infty}_{n=0}\frac{(n+1)^2}{n+2}|z|^{2n}  \quad (\mbox{by \eqref{polylog2}})\\ 
&\leq &\frac{(1-|z|^2)^4}{2}\cdot \frac{1}{(1-|z|^2)^2}\leq \frac{1}{2}.
\eeqq
This observation implies that
$$\sup_{z\in\ID}\int_{\ID}|G_z(z, \zeta)g(\zeta)|^2\, d A(\zeta)<\infty,
$$
since $g\in C(\overline{\ID})$, which is what we need.


 Finally, to complete the proof of the lemma, it remains to verify the uniqueness of $w$. If $w_1$ is also a solution to \eqref{eq-15}, then
\beqq
\begin{cases}
\Delta^2 (w-w_1)=0 & \mbox{ in $\ID$},\\
\ds w-w_1=0 & \mbox{ on $\mathbb{T}$},\\
\ds \partial_n (w-w_1)=0 & \mbox{ on $\mathbb{T}$},
\end{cases}
\eeqq
which by Theorem \Ref{pj-10} gives that $w_1=w$.
 The proof of Lemma \ref{pj-2} is complete.\epf

\noindent {\bf Proof of Theorem \ref{thm-1.1}.}
It follows from Theorems \Ref{pj-1} and \Ref{pj-10} together with Lemma \ref{pj-2} that Theorem \ref{thm-1.1} holds.
\qed

\section{Bi-Lipschitz continuity of solutions to IBDP \eqref{eq-16}}\label{sec-4}

First, we prove the Lipschitz continuity of $f$ implies the Lipschitz continuity of $u$.

\begin{lem}\label{lem4}
Let $u$ be defined by \eqref{eq-12}. Assume that $f$ satisfies the Lipschitz condition
$$|f(e^{i\theta})-f(e^{i\varphi})|\leq L |e^{i\theta}-e^{i\varphi}|
$$
and that $h\in C(\mathbb{T})$. Then for $z_1, z_2\in\ID$,
$$|u(z_1)-u(z_2)|\leq\sup_{z\in\ID}\|\nabla u(z)\|\cdot|z_1-z_2|\leq 4\Big(\frac{55}{3}L+\|h\|_{\infty,\;\mathbb{T}}\Big)|z_1-z_2|.
$$
That is, $u$ is Lipschitz continuous in $\ID$.
\end{lem}
\bpf
To prove this lemma, obviously, it suffices to show the boundedness of $\nabla u$ in $\ID$ since for $z_1, z_2\in\ID$,
\be\label{pj-16}
|u(z_{2})-u(z_{1})|=\left|\int_{[z_{1},z_{2}]}\big(u_{z}(z)\,dz+u_{\overline{z}}(z)\,d\overline{z}\big)\right|
\leq \sup_{z\in\ID}\|\nabla u(z)\|\cdot|z_1-z_2|,
\ee
where $[z_{1},z_{2}]$ stands for the segment in $\ID$ with  end points $z_{1}$ and $z_{2}$.
It follows from \eqref{pj-11} that we only need to demonstrate the boundedness of  $u_z$ and $u_{\overline{z}}$.
We first prove the boundedness of $u_z$, which is formulated in the following claim.

\bcl\label{pj-15}
For $z\in \ID$,
\be\label{pj-4}
\left|\frac{\partial u}{\partial z}(z)\right|\leq \frac{110}{3}L+2\|h\|_{\infty,\;\mathbb{T}}.
\ee
\ecl

For any $z\in \mathbb{D}$ and $\varphi\in [0, 2\pi]$,
it follows from \eqref{eq-12} and Lemma \ref{lem-2} that
\beqq
u(z) &=&\frac{1}{2\pi}\int^{2\pi}_{0}F_0(ze^{-i\theta})\big(f(e^{i\theta})-f(e^{i\varphi})\big)
\,d\theta +\frac{1}{2\pi}\int^{2\pi}_{0}H_0(ze^{-i\theta})h(e^{i\theta})\,d\theta\\
&&+f(e^{i\varphi}),
\eeqq
and then Lemma \ref{lem-4.1} implies that
$$
\frac{\partial u}{\partial z}(z)=\frac{1}{2\pi}\int^{2\pi}_{0}\frac{\partial F_0(ze^{-i\theta})}{\partial z}\big(f(e^{i\theta})-f(e^{i\varphi})\big)
\,d\theta +\frac{1}{2\pi}\int^{2\pi}_{0}\frac{\partial H_0(ze^{-i\theta})}{\partial z}h(e^{i\theta})\,d\theta.
$$

To finish the proof of \eqref{pj-4}, we let $z=re^{i\varphi}$, where $r\in[0, 1)$. Then we deduce from Lemma \ref{pj-100} together with
$1-|z|^2\leq |1-ze^{-i\theta}|(1+|z|)
$
that
\beqq\label{pj-5}
\left|\frac{\partial F_0(ze^{-i\theta})}{\partial z}\right|
\leq \frac{(1-r^2)(1+3r)}{2|1-ze^{-i\theta}|^{2}}+\frac{(1-r^2)^2(2+5r)}{2|1-ze^{-i\theta}|^{4}}
\eeqq
and
\beqq\label{pj-6}
\left|\frac{\partial H_0(ze^{-i\theta})}{\partial z}\right|\leq \frac{(1-r^2)(1+3r)}{2|1-ze^{-i\theta}|^{2}},
\eeqq
 from which we conclude that
\beqq
\left|\frac{\partial u}{\partial z}(z)\right|
&\leq&\frac{1}{2\pi}\int^{2\pi}_{0}\left|\frac{\partial F_0(z e^{-i\theta})}{\partial z}\right||f(e^{i\theta})-f(e^{i\varphi})|\,d\theta\\
&&+\frac{1}{2\pi}\int^{2\pi}_{0}\left|\frac{\partial H_0(z e^{-i\theta})}{\partial z}\right||h(e^{i\theta})|\,d\theta\\
&\leq&I_1+I_2+I_3,
\eeqq
where
$$I_1=\frac{1}{2\pi}\int^{2\pi}_{0}\frac{(1-r^2)(1+3r)}{2|1-ze^{-i\theta}|^{2}}|f(e^{i\theta})-f(e^{i\varphi})|\,d\theta,
$$
$$I_2=\frac{1}{2\pi}\int^{2\pi}_{0}\frac{(1-r^2)^2(2+5r)}{2|1-ze^{-i\theta}|^{4}}|f(e^{i\theta})-f(e^{i\varphi})|\,d\theta
$$
and
$$I_3=\frac{1}{2\pi}\int^{2\pi}_{0}\frac{(1-r^2)(1+3r)}{2|1-ze^{-i\theta}|^{2}}|h(e^{i\theta})|\,d\theta.
$$

Obviously, \eqref{pj-14}, together with the estimate $|f(e^{i\theta})-f(e^{i\varphi})|\leq 2L,$
gives
$$I_1\leq 4L\;\;\mbox{and}\;\; I_3\leq2\|h\|_{\infty,\;\mathbb{T}}.
$$

In order to estimate $I_2$, we split $\theta\in [0, 2\pi]$ into two subsets
$$E_1=\{\theta: |e^{i\theta}-e^{i\varphi}|\leq1-r\}
~\mbox{ and }~
E_2=\{\theta:\, |e^{i\theta}-e^{i\varphi}|>1-r\},
$$
where $z=re^{i\varphi}\in \ID$.  Then
$$\frac{1}{2\pi}\int^{2\pi}_{0}\frac{|e^{i\theta}-e^{i\varphi}|}{|1-z e^{-i\theta}|^{4}}\,d\theta
=\frac{1}{2\pi}\int_{E_1}\frac{|e^{i\theta}-e^{i\varphi}|}{|1-z e^{-i\theta}|^{4}}\,d\theta
+\frac{1}{2\pi}\int_{E_2}\frac{|e^{i\theta}-e^{i\varphi}|}{|1-z e^{-i\theta}|^{4}}\,d\theta.
$$
First, we find  that
$$\frac{1}{2\pi}\int_{E_1}\frac{|e^{i\theta}-e^{i\varphi}|}{|1-z e^{-i\theta}|^{4}}\,d\theta \leq \frac{1-r}{2\pi(1-r)^4} \int_{E_1}\,d\theta
\leq\frac{2\arcsin\frac{1-r}{2}}{\pi(1-r)^3}.
$$
Secondly, since $|e^{i\theta}-z|\geq1-r$ and
$$|e^{i\theta}-e^{i\varphi}|\leq|e^{i\theta}-z|+|e^{i\varphi}-z|=|e^{i\theta}-z|+1-r,
$$
we deduce that
$|e^{i\theta}-e^{i\varphi}|\leq 2|e^{i\theta}-z|,
$
and thus,
$$\frac{1}{2\pi}\int_{E_2}\frac{|e^{i\theta}-e^{i\varphi}|}{|1-z e^{-i\theta}|^{4}}\,d\theta
\leq\frac{1}{2\pi}\int^{2\pi}_{0}\frac{2d\theta}{|1-z e^{-i\theta}|^{3}} \leq \frac{2\sqrt{1+r^2}}{(1-r^2)^2},
$$
where the second inequality is a consequence of \eqref{pj-14a}.
Hence
$$I_2\leq L(1+r)^2(2+5r)\left(\frac{\arcsin\frac{1-r}{2}}{\pi(1-r)}+\frac{\sqrt{1+r^2}}{(1+r)^2}\right)\leq\frac{98}{3}L.
$$

Now, using the estimates for $I_1$, $I_2$ and $I_3$, we conclude that
$$\left|\frac{\partial u}{\partial z}(z)\right|\leq 
\frac{110}{3}L+2\|h\|_{\infty,\;\mathbb{T}},
$$
as required.
\medskip

Now, we are ready to finish the proof of the lemma based on Claim \ref{pj-15}. It follows from Lemma \ref{lem-4.1} that
 \beqq
u_{\overline{z}}(z)
&=& \frac{1}{2\pi}\int^{2\pi}_{0}\frac{\partial F_0(ze^{-i\theta})}{\partial \overline{z}}\big(f(e^{i\theta})-f(e^{i\varphi})\big)
\,d\theta +\frac{1}{2\pi}\int^{2\pi}_{0}\frac{\partial H_0(ze^{-i\theta})}{\partial \overline{z}}h(e^{i\theta})\,d\theta\\
&=& \frac{1}{2\pi}\int^{2\pi}_{0}\overline{\left(\frac{\partial F_0(ze^{-i\theta})}{\partial z}\right)}\big(f(e^{i\theta})-f(e^{i\varphi})\big)
\,d\theta \\
&& \hspace{1cm}+\frac{1}{2\pi}\int^{2\pi}_{0}\overline{\left(\frac{\partial H_0(ze^{-i\theta})}{\partial z}\right)}h(e^{i\theta})\,d\theta
\eeqq
and thus, Claim \ref{pj-15} leads to
$$\left|\frac{\partial u}{\partial \overline{z}}(z)\right|
\leq \frac{110}{3}L+2\|h\|_{\infty,\;\mathbb{T}}.
$$
Finally, we conclude the proof of the lemma from \eqref{pj-16}.
\epf

Now we concern the Lipschitz continuity of $w$.

\begin{lem}\label{lem3}
Assume that $w=G[g]$ and $g\in C(\overline{\ID})$.
Then for $z_1, z_2\in\ID$,
$$|w(z_1)-w(z_2)|\leq\sup_{z\in\ID}\|\nabla w(z)\|\cdot|z_1-z_2|\leq \frac{23}{3}\|g\|_{\infty}\cdot|z_1-z_2|.
$$
which implies the Lipschitz continuity of $w$ in $\ID$.
\end{lem}
\bpf
It follows from Lemma \ref{lem-2.6} that
\beqq\label{eq-24}
\left|w_z(z)\right|= \left|\int_{\ID}G_z(z, \zeta)g(\zeta)\,d A(\zeta)\right|\leq\frac{23}{6}\|g\|_{\infty}.
\eeqq
and
\beqq\label{eq-24a}
\left|w_{\overline{z}}(z)\right|= \left|\int_{\ID}G_{\overline{z}}(z, \zeta)g(\zeta)\, d A(\zeta)\right|\leq\frac{23}{6}\|g\|_{\infty}.
\eeqq
We conclude from \eqref{pj-11} that
$\|\nabla w(z)\|\leq (23/3)\|g\|_{\infty},
$
which completes the proof.
\epf

\noindent {\bf Proof of  Theorem \ref{thm-1.2}}\quad By Lemmas \ref{lem4} and \ref{lem3}, we obtain that
$$\|\nabla\Phi(z)\|\leq\|\nabla u(z)\|+\|\nabla w(z)\|\leq P,
$$
where $P$ is given by \eqref{eq-P}. Hence it follows from \eqref{pj-16} that for $z_1, z_2\in\ID$,
$$|\Phi(z_1)-\Phi(z_2)|\leq \sup_{z\in\ID}\|\nabla\Phi(z)\|\cdot|z_1-z_2| \leq P|z_1-z_2|.
$$
Moreover, since $A=|\Phi_{z}(0)|^2$ and $B=|\Phi_{\overline{z}}(0)|^2$, we have
$$Q=|\Phi_{z}(0)|^2-|\Phi_{\overline{z}}(0)|^2=\|\nabla\Phi(0)\|\cdot l(\nabla\Phi(0))\leq Pl(\nabla\Phi(0)),
$$
which gives that $l(\nabla\Phi(0))\geq Q/P.$ Finally, for any $z_{1}$, $z_{2}\in \ID$,  we have
\beqq
|\Phi(z_{1})-\Phi(z_{2})|&=&\left|\int_{[z_{1},z_{2}]}\left (\Phi_{z}(z)\, dz+\Phi_{\overline{z}}(z)\, d\overline{z}\right )\right|\\
&\geq&\left|\int_{[z_{1},z_{2}]}\left (\Phi_{z}(0)\, dz+\Phi_{\overline{z}}(0)\, d\overline{z}\right )\right|\\
&&-\left|\int_{[z_{1},z_{2}]}\left ( [\Phi_{z}(z)-\Phi_{z}(0)]\,dz+[\Phi_{\overline{z}}(z)-\Phi_{\overline{z}}(0)]\,d\overline{z}\right )\right|\\
&\geq& \Big(\frac{Q}{P}-2P\Big)|z_1-z_2|.
\eeqq
Hence the proof of this theorem is complete.
\qed

\subsection*{Acknowledgements}
The authors thank the referee for his/her careful reading and many useful comments.
The work of Mrs. P. Li  was supported by Centre for International Co-operation in Science (CICS)
through the award of ``INSA JRD-TATA Fellowship." The work was completed
during her visit to the Indian Statistical Institute (ISI), Chennai Centre.
The second author is on leave from the IIT Madras.


\end{document}